\documentclass[twoside]{article}
\newcommand{\bt}{\begin{theorem}}
\newcommand{\et}{\end{theorem}}
\newtheorem{lemma}{Lemma}
\newtheorem{theorem}[lemma]{Theorem}
\newtheorem{corollary}[lemma]{Corollary}

\newcommand{\tms}{\times}

\newcommand {\be}{\begin{equation}}
\newcommand {\ee}{\end{equation}}

\newcommand{\vsp}{\vskip 1em}

\def \qed {\hfill \vrule height6pt width 6pt depth 0pt}
\def \qed {\hfill \vrule height6pt width6pt depth0pt}

\def \textit{\it}
\def \bt{\begin{theorem}}
\def \et{\end{theorem}}
\def \bl{\begin{lemma}}
\def \el{\end{lemma}}
\def \bc{\begin{corollary}}
\def \ec{\end{corollary}}
\def \be{\begin{equation}}
\def \ee{\end{equation}}

\def \text{\mbox}

\def \noi{\noindent}

\title{Monotonicity properties of certain Laplacian eigenvectors associated with trees}
\author{ Ravindra B. Bapat\footnote{
The  author acknowledges support from the JC Bose Fellowship, Department
of Science and Technology, Government of India.}
  \\ [.1cm]
Indian Statistical Institute,\\ Delhi Centre, 7 S.J.S.S. Marg, \\  New Delhi 110 016, India\\
 {rbb@isid.ac.in}}

\date{}

\begin{document}
\maketitle

\begin{abstract}
Nath and Paul (Linear Algebra Appl., 460 (2014), 97-110)
have shown that the largest distance Laplacian
eigenvalue of a path is simple and the corresponding eigenvector has
properties similar to the Fiedler vector. We given an alternative proof,
establishing a more general result in the process. It is conjectured that
a similar phenomenon holds for any tree.
\end{abstract}

\vspace{3mm}

\section{Introduction}

Let $T$ be a tree with vertex set $V(T) = \{1, \ldots, n\}.$
We denote the degree of vertex $i$ by $\delta_i, i = 1, \ldots, n.$
         Recall that the Laplacian matrix  $L$ of $T$ is the $n \tms n$ matrix with its $(i,j)$-element, $i \not = j,$
 equal to $-1,$ if $i$ and $j$ are adjacent, and zero otherwise. The diagonal elements of $L$  are $\delta_1, \ldots, \delta_n.$
 It is well-known that
 $L$ is a positive semidefinite matrix with rank $n-1.$ Let
 $0 = \lambda_1 <  \lambda_2 \le \cdots \le  \lambda_n$ be the eigenvalues of $L.$
 The second smallest eigenvalue $\lambda_2$ is termed the {\it algebraic
 connectivity} of $T$ and the corresponding eigenvector is a Fiedler vector.
 For basic properties of the Laplacian matrix we refer to \cite{bapat}.
 We state the following classical result (see, for example,
 \cite{bapat}, Chapter 8, \cite{jason}, Chapter 6).

 \bt \label{fiedler} Let $T$ be a tree with $V(T) = \{1, \ldots, n\}$ and let $f$ be a Fiedler
 vector of $T.$ Let $f(v)$ denote the component of $f$ indexed by the vertex $v.$
 Then one of the two cases must occur:

 \begin{description}
 \item {Case (i).} $f(v) \not = 0$ for any $v.$ Then there is a unique edge
   $\{u_1,u_2\}$ such that $f(u_1) > 0$ and $f(u_2) < 0.$ Moreover, the values of $f$ increase along any path starting at $u_1$ and not containing $u_2,$ and the values of $f$ decrease along any path starting at $u_2$ and not containing $u_1.$

 \item {Case (ii).} $f(v) = 0$ for some v. Then there is a unique vertex $u$
  such that $f(u) = 0$ and $u$ is adjacent to a vertex
 on which $f$ takes a nonzero value. Moreover, along any path starting at $u,$
 the values of $f$ either increase, decrease or are identically zero.

 \end{description}
 \et

The unique edge in Case (i) is called the {\it characteristic edge}, while the unique vertex in Case (ii) is called the {\it characteristic vertex.}

Let $T$ be a tree with $V(T) = \{1, \ldots, n\}.$
The distance between vertices $i,j \in V(T),$ denoted $d_{ij},$ is defined to
be the length (the number of edges) of the (unique) path from $i$ to $j.$
 We set $d_{ii} = 0, i = 1, \ldots, n.$
The distance matrix $D$ is the $n \tms n$ matrix with $(i,j)$-element
equal to $d_{ij}.$

 The distance Laplacian of $T,$ denoted $D^L,$
is the $n \tms n$ matrix with $(i,j)$-element $-d_{ij}$ if $i \not = j,$
and the $(i,i)$-element equal to $\sum_{j=1}^n d_{ij}, i = 1, \ldots, n.$
It has been conjectured by Nath and Paul \cite{nathpaul} that if $f$ is an
eigenvector corresponding to the largest eigenvalue of $D^L,$ then it has properties
similar to that of a Fiedler vector, more specifically, it satisfies
either Case (i) or Case (ii)
of Theorem \ref{fiedler}. The conjecture was confirmed for a path in \cite{nathpaul}.
In this note we prove a more general statement for a path, and conjecture that
a similar phenomenon occurs for any tree.

\section{Eigenvectors of certain Laplacians associated with trees}

Let $A$ be a symmetric $n \tms n$ matrix. The Laplacian associated with $A,$
denoted $A^L,$ is the $n \tms n$ matrix with $(i,j)$-element $-a_{ij}$ if $i \not = j,$
and the $(i,i)$-element equal to $\sum_{j \not = i} a_{ij}, i = 1, \ldots, n.$
We formulate the following conjectures, based on extensive numerical
computations.

\vsp

{\bf Conjecture 1}
Let $T$ be a tree with $V(T) = \{1, \ldots, n\}.$ Let $A$ be a symmetric, nonnegative
 $n \tms n$ matrix such that for any distinct $i,j$ and $k,$
$a_{ij} \le a_{ik},$ whenever $j$ and $k$ are adjacent and $d_{ij} > d_{ik}.$
Let $0 = \lambda_1 \le  \lambda_2 \le \cdots \le  \lambda_n$ be the eigenvalues of $A^L.$ Then there exists an eigenvector $f$  of $A^L$ corresponding to $\lambda_2$
such that $f$ satisfies either Case (i) or Case (ii) of Theorem \ref{fiedler}.

\vsp

{\bf Conjecture 2}
Let $T$ be a tree with $V(T) = \{1, \ldots, n\}.$ Let $A$ be a symmetric, nonnegative $n \tms n$ matrix such that for any distinct $i,j$ and $k,$
$a_{ij} \ge a_{ik},$ whenever $j$ and $k$ are adjacent and $d_{ij} > d_{ik}.$
Then there exists an eigenvector $f$ of $A^L$ corresponding to its largest eigenvalue
such that $f$ satisfies either Case (i) or Case (ii) of Theorem \ref{fiedler}.

\vsp

The motivation for these conjectures is as follows. If $A$ is the adjacency matrix
of the tree, than it satisfies the condition in Conjecture 1 and hence
either Case (i) or Case (ii) are satisfied by Theorem \ref{fiedler}.
If $A$ is the distance matrix of the tree, then
it satisfies the condition in Conjecture 2 and we obtain the conjecture
of Nath and Paul \cite{nathpaul}.

In the rest of this paper we show that Conjectures 1 and 2
are true when the tree $T$ is a path.
The proof of Conjecture 1 for a path is essentially contained in
\cite{boman}, Theorem 3.2. We give the proof for completeness.
We use the same proof technique and prove Conjecture 2 for a path.
Let ${\bf 1}$ denote the vector of appropriate size of all ones.
We now state a preliminary result from \cite{boman}.
The proof is easy and is omitted.

\bl \label{lem7} Let $S$ and $T$ be matrices of order  $(n-1) \tms n$ and $n \tms (n-1)$
respectively, given by
$$S = \left(
        \begin{array}{rrrrr}
          -1 & 1 & 0 & 0 & 0 \\
          0 & -1 & 1 & 0 & 0 \\
          \vdots &  &  &  & \vdots \\
          0 & 0 & \cdots & -1 & 1 \\
        \end{array}
      \right),\,\,
 T = \left(
       \begin{array}{cccc}
         0 & 0 & \cdots & 0 \\
         1 & 0 & \cdots & 0 \\
         1 & 1 & \cdots & 0 \\
         \vdots &  &  & \vdots \\
         1 & 1 & \cdots & 1 \\
       \end{array}
     \right).$$
     Then $ST = I_{n-1}$ and $TS = I_n - {\bf 1}{e_1}',$
     where $e_1$ is the first column of the identity matrix $I_n.$
     \el

We now prove a reformulation of a result from \cite{boman}.

 \bl \label{lem8} Let $A$ be an $n \tms n$ symmetric matrix and
 let $M = SA^LT,$ where $S$ and $T$ are as in Lemma \ref{lem7}.
Let $0 = \lambda_1,  \lambda_2,  \cdots,  \lambda_n$ be the eigenvalues of $A^L.$ Then

\begin{description}

\item {(i)} The eigenvalues of $M$ are $\lambda_2, \ldots, \lambda_n.$

\item {(ii)} If $y$ is an eigenvector of $A^L$ such that $y \perp {\bf 1},$
then $Sy$ is an eigenvector of $M$ corresponding to the same eigenvalue.

\item {(iii)} If $y$ is an eigenvector of $M,$ then there exists an eigenvector
$x$ of $A^L$ corresponding to the same eigenvalue such that $x \perp {\bf 1}$ and $y = Sx.$

\end{description}

\el

\noi
{\bf Proof:}
(i).  Recall that if $X$ and $Y$ are matrices of order $p \tms q$
and $q \tms p$ respectively, where $p \le q,$ and if $\mu_1,
\ldots, \mu_p$ are the eigenvalues of $XY,$ then the eigenvalues of $YX$ are
given by $\mu_1, \ldots, \mu_p,$ along with $0$ repeated $q-p$ times.
Note that $M = SA^LT$ and
$A^LTS = A^L(I_n - {\bf 1}{e_1}') =  A^L,$ since $A^L{\bf 1} = 0.$
It follows that     the eigenvalues of $M$ are $\lambda_2, \ldots, \lambda_n.$

(ii). Let $A^Ly = \alpha y.$ Note that since $y \not = 0$
and $y \perp {\bf 1},$ then $Sy \not = 0.$ We have
  $MSy = SA^LTSy = SA^L(I_n - {\bf 1}{e_1}')y = SA^Ly = \alpha Sy$
  and the proof is complete.

  (iii). If $My = \alpha y,$ then letting $x = Ty,$ we have
$  SA^LTSx = \alpha Sx,$ which implies $SA^Lx = \alpha Sx.$
Since $x \perp {\bf 1},$ we conclude $A^L x = \alpha x.$
  \qed

\bl \label{lem9} Let $A$ be a symmetric $n \tms n$ matrix and let $M = SA^LT,$ where
$S$ and $T$ are as in Lemma \ref{lem7}.
\begin{description}

\item {(i)} Suppose $a_{ij} \ge a_{ik}$ if $j < k < i$ and $a_{ij} \le a_{ik}$ if $i < j < k.$ Then $m_{ij} \ge 0$ for any $i \not = j.$

\item {(ii)} Suppose $a_{ij} \le a_{ik}$ if $j < k < i$ and $a_{ij} \ge a_{ik}$ if $i < j < k.$ Then $m_{ij} \le 0$ for any $i \not = j.$

 \end{description}
 \el

  \noi
  {\bf Proof:}

  (i).     Suppose $a_{ij} \ge a_{ik}$ if $j < k < i$ and $a_{ij} \le a_{ik}$ if $i < j < k.$ For $i < j,$
  \begin{eqnarray*}
  m_{ij} &=& (SA^LT)_{ij} \\
  &=& \sum_{k=1}^n (SA^L)_{ik}t_{kj} \\
    &=& \sum_{k=j+1}^n (a_{ik} - a_{i+1,k}) \\
  &\ge& 0.
  \end{eqnarray*}

For $i > j,$ using the fact that $A^L$ has zero row sums,
  \begin{eqnarray*}
  m_{ij} &=&  \sum_{k=j+1}^n (a_{ik} - a_{i+1,k}) \\
  &=& \sum_{k=1}^j (-a_{ik} + a_{i+1,k}) \\
    &\ge& 0.
  \end{eqnarray*}

The proof of (ii) is similar. \qed

An examination of the proof of Lemma \ref{lem9} reveals that the following
result is true, which we state without proof.

\bl \label{lem10} Let $A$ be a symmetric $n \tms n$ matrix and let $M = SA^LT,$ where
$S$ and $T$ are as in Lemma \ref{lem7}.
\begin{description}

\item {(i)} Suppose $a_{ij} > a_{ik}$ if $j < k < i$ and $a_{ij} < a_{ik}$ if $i < j < k.$ Then $m_{ij} > 0$ for any $i \not = j.$

\item {(ii)} Suppose $a_{ij} < a_{ik}$ if $j < k < i$ and $a_{ij} > a_{ik}$ if $i < j < k.$ Then $m_{ij} < 0$ for any $i \not = j.$

 \end{description}
 \el

\bl \label{lem11}
Let $A$ be a symmetric $n \tms n$ matrix such that $a_{ij} \ge 0$ for all
$i \not = j.$ If $\lambda_n$ is the largest eigenvalue of $A,$ then $A$
has an eigenvector $x$ for $\lambda_n$ with $x_i \ge 0, i = 1,
\ldots, n.$ Furthermore, if $a_{ij} > 0$ for all $i \not = j,$
then $x_i > 0, i = 1, \ldots, n.$
\el

{\bf Proof:}
The result follows from the Perron-Frobenius Theorem,
applied to the matrix $A + \beta I$ for a sufficiently
large $\beta.$ \qed

\vsp

The following is our main result.

\bt \label{thm11}
Let $A$ be a symmetric, nonnegative $n \tms n$ matrix with $a_{ij} \ge 0$
for all $i \not = j$ and $a_{ii} = 0, i = 1, \ldots, n.$
 Let $0 = \lambda_1 \le  \lambda_2 \le \cdots \le  \lambda_n$ be the eigenvalues of $A^L.$

\begin{description}
\item {(i)} Suppose $a_{ij} \ge a_{ik}$ if $j < k < i$ and $a_{ij} \le a_{ik}$
if $i < j < k.$ Then $A^L$ has an eigenvector $x$ corresponding to $\lambda_n$
such that $x \perp {\bf 1}$  and $x_1 \le \cdots \le x_n.$

\item {(ii)} Suppose $a_{ij} \le a_{ik}$ if $j < k < i$ and $a_{ij} \ge a_{ik}$
if $i < j < k.$ Then $A^L$ has an eigenvector $x$ corresponding to $\lambda_2$
such that $x \perp {\bf 1}$ and $x_1 \le \cdots \le x_n.$

\end{description}

\et

\noi
{\bf Proof:} (i). Let $M = SA^LT,$ where $S$ and $T$ are as in Lemma \ref{lem7}.
By Lemma \ref{lem9}, $m_{ij} \ge 0$ for all $i \not = j.$ By Lemma \ref{lem8},
the largest eigenvalue of $M$ is $\lambda_n,$ and by Lemma \ref{lem11}, $M$
has an eigenvector $y$ corresponding to $\lambda_n$ with $y_i \ge 0,
i = 1, \ldots, n.$ By Lemma \ref{lem8}, $A^L$ has an eigenvector $x \perp {\bf 1}$
corresponding to $\lambda_n$ such that $y = Sx.$ Since $y_i \ge 0,
i = 1, \ldots, n,$ it follows that $x_1 \le x_2 < \cdots  \le x_n.$

(ii). Let $M$ be as in (i). By Lemma \ref{lem9}, $m_{ij} \le 0$ for all $i \not = j.$
Let $\beta$ be sufficiently large so that $\beta I - M$ has all entries
nonnegative. By Lemma \ref{lem8}, the eigenvalues of $M$ are $\lambda_2,
\ldots, \lambda_n,$ and hence $\beta I - M$ has eigenvalues $\beta - \lambda_2,
\ldots, \beta - \lambda_n.$ The largest eigenvalue of $\beta I - M$
is $\beta -\lambda_2$ and by the Perron-Frobenius Theorem it has an eigenvector
$y$ corresponding to this eigenvalue with $y_i \ge 0, i = 1, \ldots, n.$
Clearly $y$ must also be an eigenvector of $M$ corresponding to $\lambda_2.$
By Lemma \ref{lem8}, $A^L$ has an eigenvector $x \perp {\bf 1}$
corresponding to $\lambda_2$ such that $y = Sx.$ Since $y_i \ge 0,
i = 1, \ldots, n,$ it follows that $x_1 \le x_2 < \cdots  \le x_n.$
\qed

\bc \cite{nathpaul} Let $T$ be the path with $V(T) = \{1, \ldots, n\}$
and $E(T) = \{\{i,i+1\}: 1 \le i \le n-1\}.$ Let $D$ be the distance matrix of
$T$ and let $\lambda_n$ be the largest eigenvalue of $D^L$ with a corresponding
eigenvector $x.$ Then either $x_1 \ge \cdots \ge x_n$ or $x_1 \le \cdots \le x_n.$
\ec

\noi
{\bf Proof:} Note that $D$ satisfies the condition in Theorem \ref{thm11}, (i).
           By Lemma \ref{lem10}, if $M = SD^LT,$ then $m_{ij} > 0$ for all $i \not = j$
and hence $\lambda_n$ is a simple eigenvalue of $M,$ and therefore of $D^L.$
It follows that the eigenvector of $D^L$ corresponding to $\lambda_n$
must be unique up to a scalar multiple. By Theorem \ref{thm11}, for any eigenvector $x$
corresponding to $\lambda_n,$  either $x_1 \ge \cdots \ge x_n$ or $x_1 \le \cdots \le x_n.$
\qed


\begin{thebibliography}{MM}

\bibitem{boman} J.E. Atkins, E.G. Boman and B. Hendrickson, A spectral algorithm
for seriation and the consecutive ones problem, {\em SIAM J. Comput.}
{\bf 28(1)} (1998), 297--310.


\bibitem{bapat} R.B. Bapat, {\em Graphs and matrices}, Second Edition,
Springer, London, 2014.

\bibitem{jason} Jason J. Molitierno, {\it Applications of combinatorial matrix theory
to Laplacian matrices of graphs}, CRC Press, 2012.

\bibitem{nathpaul} Milan Nath and Somnath Paul, On the distance Laplacian spectra
of graphs, {\em Linear Algebra Appl.} {\bf 460} (2014), 97--110.


\end{thebibliography}
\end{document}